\documentclass[final]{IEEEtran}

\usepackage{amsmath,amsfonts,amssymb,amsopn}
\usepackage{graphicx}
\usepackage{xcolor}
\usepackage{cite}
\usepackage{xspace}
\usepackage[caption=false,font=footnotesize]{subfig}
\usepackage{mathtools}
\usepackage{url}
\usepackage{booktabs}
\usepackage{algorithm}
\usepackage{algpseudocode}
\usepackage{ifpdf}
\ifpdf
  \usepackage[pdftex,colorlinks=true,breaklinks=true,citecolor=blue]{hyperref}
\else
  \usepackage[ps2pdf,colorlinks=true,citecolor=blue]{hyperref}
\fi
\usepackage{float}
\floatname{algorithm}{Procedure}
\def\NoNumber#1{{\def\alglinenumber##1{}\State #1}\addtocounter{ALG@line}{-1}}

\ifCLASSOPTIONonecolumn
    
    \renewcommand{\baselinestretch}{2}
    \usepackage[margin=25mm]{geometry}
\else
    
\fi

\newcommand{\bv}[1]{\boldsymbol{#1}}
\newcommand{\dfn}{\triangleq}
\newcommand{\untsph}{\mathbb{S}^{2}}
\newcommand{\shc}[3]{({#1})_{#2}^{#3}}
\newcommand{\lsph}{L^2(\untsph)}
\newcommand{\conj}[1]{\overline{#1}}
\newcommand{\unit}[1]{\bv{\hat{#1}}}
\newcommand{\lsphL}[1]{\mathcal{H}_{#1}}
\newcommand{\intsph}{\int_{\untsph}}
\newcommand{\inttheta}{\int_{\theta=0}^{\pi}}
\newcommand{\intphi}{\int_{0}^{2\pi}}
\newcommand{\figref}[1]{Fig.\,\ref{#1}}
\newcommand{\secref}[1]{Section\,\ref{#1}}
\newcommand{\appref}[1]{Appendix\,\ref{#1}}
\newcommand{\tabref}[1]{Table\,\ref{#1}}
\newcommand{\tilf}{\tilde f}

\newcommand{\ndh}{N_{\textrm{DH}}}
\newcommand{\nmw}{N_{\textrm{MW}}}
\newcommand{\ngl}{N_{\textrm{GL}}}
\newcommand{\plms}{\tilde{P}}
\newcommand{\nps}{\ensuremath{N_\textrm{O}}}
\DeclarePairedDelimiterX\innerp[2]{\langle}{\rangle}{#1,#2}

\newcommand{\sfn}[2]{\Lambda^{#1}_{#2}}          
\newcommand{\amp}{\Lambda}                        
\newcommand{\kap}{\kappa}
\newcommand{\Pm}{\mathbf{P}_m}
\newcommand{\Psm}{\mathbf{P}_m^{\scriptscriptstyle <}}   
\newcommand{\Lm}{\bv{\Lambda}_m}
\newcommand{\Vm}{\mathcal{V}_{m}}

\newcommand{\abs}[1]{\lvert #1 \rvert}
\newcommand{\norm}[1]{\lVert #1 \rVert}

\newcommand{\placeholder}[2]{\fbox{\parbox[c][#1][c]{0.96\linewidth}{\centering\small\textit{Figure placeholder:} #2}}}

\newtheorem{lemma}{Lemma}

\newtheorem{remark}{Remark}
\newtheorem{theorem}{Theorem}
\graphicspath{{figs/},{pdfs/},{pdf/},{Figures/}}
\begin{document}
\title{An $O(L^3)$ Spherical Harmonic Transform on $L^2$ Samples}

\author{%
Zubair~Khalid,~\IEEEmembership{Senior Member,~IEEE}%
\thanks{Z.~Khalid is with the Department of Electrical Engineering, Syed Babar Ali School of Science and Engineering, Lahore University of Management Sciences, Lahore 54792, Pakistan (e-mail: zubair.khalid@lums.edu.pk).}
\thanks{An open-source implementation of the transforms of this paper, the optimal-dimensionality SHT package (\texttt{osht}), is publicly available at \url{https://github.com/zubairkhalid1/osht/}.}
}
\maketitle

\begin{abstract}
For the representation of signals on the sphere band-limited at $L$ using the optimal number of samples, $L^2$, equal to the number of degrees of freedom in harmonic space and placed on $L$ iso-latitude rings, we develop a spherical harmonic transform (SHT) that is theoretically exact, requires no solution of linear systems and no precomputation beyond the placement of the rings, and has complexity $O(L^3)$. The transform rests on the observation that, for each order $m$, the samples on the rings that resolve order $m$ are a complete set of samples for the order-$m$ profile of the signal, so that its values on the remaining rings follow from a closed-form Lagrange basis, the sampling functions of order $m$. We prove that the proposed optimal-dimensionality SHT (OSHT) recovers every band-limited signal exactly in exact arithmetic, so that the $L^2$ samples admit a sampling theorem in the collocation sense. For the final harmonic analysis we give two exact back-ends: Gauss-Legendre quadrature of the interpolated profiles and the algorithm of the equiangular sampling theorem, which applies without change because the rings coincide with its grid. We identify the amplification constant that governs numerical error, show that it replaces the condition number used previously for the placement of rings, propose a placement that minimizes it, and give an error propagation model that explains the growth of error with the band-limit. A multi-pass refinement converges in two passes to the rounding floor of the samples. We also describe a stable log-domain evaluation of the sampling functions, a rescaled Legendre recurrence, and a forbidden-region skipping rule that halves the cost of the elimination. Numerical experiments up to $L=2048$ characterize the accuracy, measured by the maximum and mean coefficient errors and the signal-to-noise ratio, the stability and the computation time.
\end{abstract}


\section{Introduction}
\IEEEPARstart{S}{ignals} defined on the sphere arise in geodesy~\cite{Wieczorek:2007}, cosmology~\cite{Spergel:2007}, computer graphics~\cite{Ng:2004}, medical imaging~\cite{Chung:2007}, astrophysics~\cite{Jarosik:2010}, quantum chemistry~\cite{Ritchie:1999}, acoustics~\cite{Zhang:2012} and planetary science~\cite{Audet:2011}, among other fields. The spherical harmonic transform~(SHT), the spherical counterpart of the Fourier transform, maps a signal between the spatial and the harmonic~(spectral) domains and is the basic tool of signal processing on the sphere~\cite{Driscoll:1994,McEwen:2011,Kennedy-book:2013}. A sampling scheme and its associated SHT are judged by three criteria: the number of samples needed to represent a signal band-limited at $L$, the computational complexity of the transform, and its numerical accuracy~\cite{Khalid:2014}. A band-limited signal has $L^2$ degrees of freedom in harmonic space, and a sampling scheme that attains this number has optimal spatial dimensionality.

The optimal-dimensionality sampling scheme~\cite{Khalid:2014} is an iso-latitude scheme with exactly $\nps=L^2$ samples, placed on $L$ rings carrying $1,3,5,\hdots,2L-1$ equiangular points, together with a forward SHT that proceeds order by order from $|m|=L-1$ down to $0$. At each order the harmonic coefficients are obtained by solving a linear system formed on the rings that carry at least $2|m|+1$ points, and their contribution is then removed from the remaining rings, so that these rings become band-limited along longitude for the next order. Its complexity is $O(L^4)$ in theory because of the linear system solved at every order, and its accuracy depends on the conditioning of those systems, which was controlled by a greedy placement of the rings that minimizes the condition number. Subsequent work refined the placement and added a multi-pass correction~\cite{Nafees:2017}, and extended the scheme to spin functions~\cite{Elahi:2018}, to antipodally symmetric signals and multiple shells in diffusion MRI~\cite{Bates:2016,Bates:2017} and to head-related transfer function measurement~\cite{Bates:2015}.

\begin{figure*}[t]
\centering
\includegraphics[width=\textwidth]{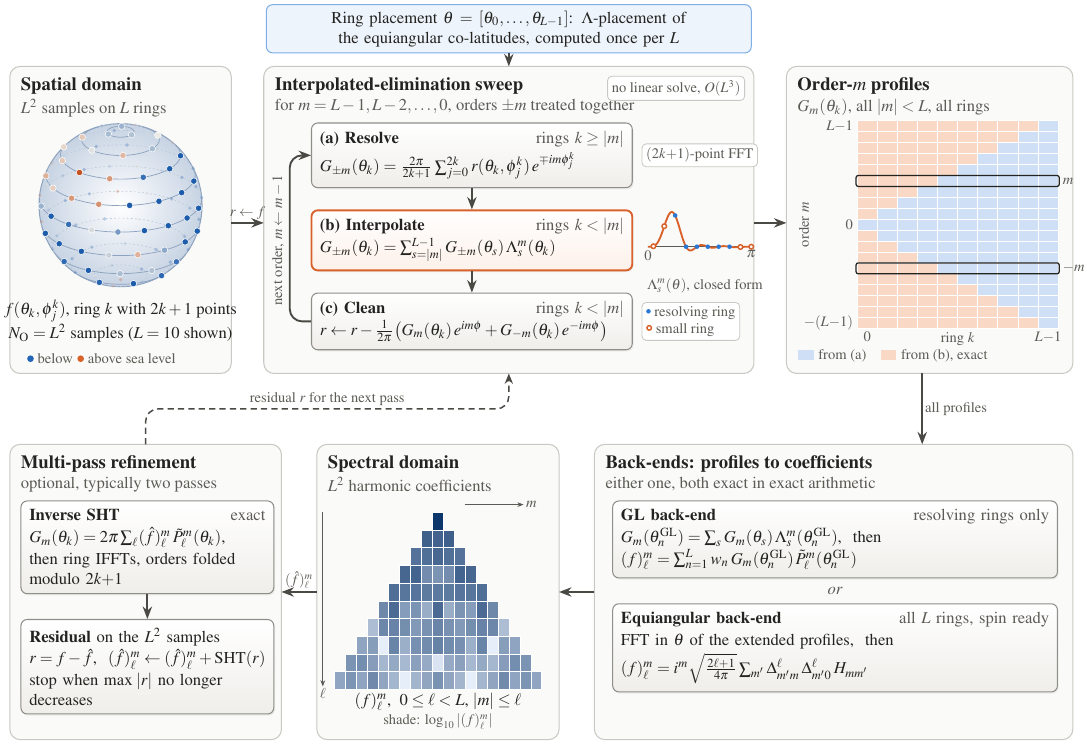}
\caption{Processing pipeline of the proposed forward OSHT of \secref{sec:transform}, illustrated for $L=10$. Top left: the $\nps=L^2$ samples $f(\theta_k,\phi^k_j)$, with $2k+1$ equiangular points on ring $k$ at the co-latitudes $\bv\theta$ of the $\amp$-placement (\secref{sec:placement}); samples on the far side of the sphere are shown faded. Top center: the sweep of Procedure~\ref{Procedure:IE} over the orders $m=L-1,\hdots,0$, starting from the residual $r=f$; at each order, step (a) resolves $G_{\pm m}$ on the rings $k\ge|m|$ by FFTs, step (b) carries these values to the rings $k<|m|$ through the sampling functions of order $m$, shown in the inset with their values at the resolving rings (filled markers) and at the small rings (open markers), and step (c) removes the orders $\pm m$ from the small rings. Top right: the profiles $G_m(\theta_k)$ after the sweep for every $|m|<L$ and every ring, obtained in step (a) (blue) or in step (b) (orange), with the rows of the current orders $\pm m$ outlined. Bottom right: the GL back-end uses the resolving rings only and the equiangular back-end all rings (Procedure~\ref{Procedure:backend}), both exact in exact arithmetic (Theorem~\ref{thm:exact}). Bottom center: the $L^2$ coefficients $\shc{f}{\ell}{m}$ of an example signal, shaded by $\log_{10}|\shc{f}{\ell}{m}|$. Bottom left: the optional multi-pass refinement of \secref{sec:multipass}, which synthesizes the current estimate by the inverse transform, forms the residual on the samples and returns it to the sweep (dashed arrow). No linear system is solved, and each pass costs $O(L^3)$.}
\label{fig:pipeline}
\end{figure*}

In this paper we develop a different transform on the same $L^2$ samples, which we refer to as the optimal-dimensionality spherical harmonic transform (OSHT). The proposed OSHT keeps the descending sweep over orders but never solves a linear system: for each order the profile of the signal on the rings that cannot resolve that order is obtained from its values on the rings that can, through a closed-form Lagrange basis that we call the sampling functions of order $m$. The resulting transform is theoretically exact, has complexity $O(L^3)$ with no precomputation, and admits a complete error analysis in terms of a single amplification constant per order. The processing pipeline, from the samples to the harmonic coefficients, is summarized in \figref{fig:pipeline}.

\subsection{Relation to Prior Work}
Iso-latitude sampling schemes~\cite{Driscoll:1994,Crittenden:1998,Healy:2003,Gorski:2005,Wiaux:2005:2,Kostelec:2008,Huffenberger:2010,McEwen:2011,Sneeuw:1994,Blais:2006} enable a separation of variables in the computation of the SHT. Among them, the sampling theorems of Driscoll and Healy~\cite{Driscoll:1994} with $\ndh\sim4L^2$ samples, of Gauss--Legendre (GL) quadrature~\cite{Skukowsky:1986,Doroshkevich:2005} with $\ngl\sim2L^2$ samples and of McEwen and Wiaux~\cite{McEwen:2011}, referred to as the equiangular sampling theorem, with $\nmw=(L-1)(2L-1)+1\sim2L^2$ samples lead to theoretically exact transforms of complexity $O(L^3)$ that require no precomputation. All three obtain the harmonic coefficients by an exact quadrature rule, which requires the integration of products of two band-limited functions, of degree $2L-2$, and therefore at least about $2L^2$ samples. The optimal-dimensionality scheme~\cite{Khalid:2014} reaches $L^2$ samples by abandoning quadrature in favor of collocation, that is, by solving for the coefficients, and this was the reason it was presented as accurate rather than exact.

The present work shows that collocation on the optimal-dimensionality scheme is exact as well, and that the collocation step can be written in closed form. The key object is the Lagrange basis of the $L^2$-dimensional space of band-limited signals on an $L^2$-point set, which was studied in~\cite{Giacaglia:1972} under the name of sampling functions as an alternative to spherical harmonics for the representation of the geoid and the geopotential. That report contains, in a notation adapted to geophysics, the point set that we use (up to a reflection through the equator), the structure of the sampling functions order by order, and the closed-form inverse of the per-order system; it did not contain a transform, a stability analysis, or the ring placement that makes the scheme usable at large band-limits. We identify this connection, and we relate the condition for orthogonal sampling functions given in~\cite{Giacaglia:1972} to tight spherical designs~\cite{An:2010,Bannai:2009,Delsarte:1977}, which shows that no $L^2$-point set with orthogonal sampling functions exists beyond the regular tetrahedron, a question left open in~\cite{Giacaglia:1972}.

The computation of SHTs has developed along four lines in the last three to four decades. The first line consists of high-performance libraries for the schemes with about $2L^2$ and $4L^2$ samples. These include libsharp~\cite{Reinecke:2013}, SHTns~\cite{Schaeffer:2013} for pseudo-spectral simulation, the butterfly-compressed Wavemoth~\cite{Seljebotn:2012} and, most recently, S2FFT~\cite{Price:2024}, which implements the equiangular and HEALPix transforms in JAX so that they are differentiable and run on GPUs. These codes reduce the constants of the $O(L^3)$ scaling for the GL, equiangular and HEALPix grids, but they do not change the number of samples. The second line consists of fast algorithms that compress the associated Legendre transform by the butterfly scheme~\cite{Tygert:2010,Bremer:2021}, or that convert between spherical harmonic expansions and bivariate Fourier series~\cite{Slevinsky:2019}. These algorithms reach a cost of $L^2$ times a small power of $\log L$, and the conversion of~\cite{Slevinsky:2019} is backward stable. The conversion to bivariate Fourier series exploits the torus structure of the sphere, as the equiangular sampling theorem~\cite{McEwen:2011} does, and it is relevant here because the equiangular back-end of \secref{sec:backends} is a transform of this kind, applied to the interpolated profiles. The third line consists of spherical designs and extremal point sets with a number of points of the order of the harmonic degrees of freedom: designs computed by optimization~\cite{Graf:2011}, well-conditioned and efficient designs with good geometric properties~\cite{An:2010,Womersley:2018}, and the extremal systems of Sloan and Womersley~\cite{Sloan:2004}. These point sets give exact quadrature or interpolation on $O(L^2)$ unstructured points, at the cost of $O(L^6)$ transforms, by the solution of a dense system of size $L^2$, and of a nonlinear search for the points, whereas the ring structure used here reduces the cost to $O(L^3)$. The fourth line concerns the optimal-dimensionality scheme itself, with the refinements and extensions recalled above~\cite{Nafees:2017,Elahi:2018,Bates:2016,Bates:2017,Bates:2015}. Iso-latitude sampling theorems also exist on the rotation group~\cite{McEwen:2015}.

In applications, taking samples on the sphere and computing the SHT from them remains the standard approach in planetary geophysics, where the SHTools package~\cite{Wieczorek:2018} analyzes gravity and topography on Driscoll--Healy and GL grids, in cosmology, where the Planck legacy maps are analyzed on the HEALPix grid~\cite{Planck:2020}, and in spatial audio, where spherical microphone arrays sample the sound field on the sphere and process it in the spherical harmonic domain~\cite{Rafaely:2019}. In each of these fields the number of samples is a cost, in instrument time, in pixels or in microphones, and a scheme that attains the lower bound of $L^2$ samples with an accurate transform is therefore of direct use. SHTs have also become a building block of learning on the sphere, from spherical convolutional networks~\cite{Cohen:2018,Esteves:2018} to spherical Fourier neural operators for weather and climate emulation~\cite{Bonev:2023}. In these networks the transform is evaluated at every layer, so that both its number of samples and its differentiability matter. The transform of this paper is a composition of FFTs, dense matrix-vector products and quadrature sums, and it is therefore differentiable in the same sense as the transforms of~\cite{Price:2024}.

The stable evaluation of Lagrange bases through barycentric formulas is classical~\cite{Berrut:2004,Higham:2004}; here it is combined with the weight $(\sin\theta)^{|m|}$ of the order-$m$ profile, and keeping the weight inside a log-domain evaluation makes the closed-form inverse comparable in accuracy to an LU factorization.

\subsection{Contributions and Paper Organization}
In broad brush terms, the contributions of the paper are as follows:
\begin{itemize}
\item We develop an $O(L^3)$ forward SHT on $L^2$ samples that solves no linear system and requires no precomputation beyond the ring co-latitudes, by replacing every per-order solve of~\cite{Khalid:2014} with the closed-form sampling functions of order $m$. We prove that the transform recovers every band-limited signal exactly in exact arithmetic, so that the optimal-dimensionality samples admit a sampling theorem in the collocation sense (Theorem~\ref{thm:exact}), and we provide two exact back-ends for the final harmonic analysis: GL quadrature of the interpolated profiles, and the equiangular transform of~\cite{McEwen:2011}.
\item We use the closed form of the per-order inverse as a weighted Lagrange basis, given in~\cite{Giacaglia:1972}, and relate the scheme to the sampling functions of~\cite{Giacaglia:1972}, to the reproducing-kernel form and to the tight-design condition. We show that the 1-3-5 point distribution of~\cite{Giacaglia:1972} is equivalent to the alternating ordering of~\cite{Khalid:2014}. We also give a stable log-domain evaluation of the sampling functions with nested barycentric weights updated in $O(L-|m|)$ per order, and we identify the unweighted route as the unstable one.
\item We identify the amplification constant $\amp_m$ that governs the error of the elimination, show that it replaces the condition number $\kap_m$ of~\cite{Khalid:2014}, bound the contribution of the back-end, and give an error propagation model that explains the growth of the single-pass error with the band-limit. On this basis we propose a placement of the rings that minimizes $\amp_m$, which lowers the error of both the new and the previous transform, and we give an $O(L^4)$ implementation of the condition-number placement of~\cite{Khalid:2014} by rank-one spectral updates.
\item We integrate the multi-pass refinement of~\cite{Nafees:2017} with the new transform and characterize its convergence. We also describe and quantify several numerical enhancements: FFT-based ring transforms, a rescaled Legendre recurrence immune to seed underflow, and a forbidden-region skipping rule that reduces the cost of the elimination by about half without affecting accuracy.
\item We carry out experiments up to $L=2048$ that characterize the single-pass and multi-pass accuracy of the transform with the equiangular back-end for both placements, measured by the maximum and mean errors and the signal-to-noise ratio, as well as its computation time.
\end{itemize}
The remainder of the paper is organized as follows. \secref{sec:models} reviews harmonic analysis on the sphere and the optimal-dimensionality sampling scheme in the notation of~\cite{Khalid:2014}. The sampling functions are developed in \secref{sec:sampling_functions}, the transform and its exactness in \secref{sec:transform}, its numerical realization and computational complexity in \secref{sec:numerics}, and the error analysis, the placement of rings and the multi-pass refinement in \secref{sec:error}. The numerical experiments are reported in \secref{sec:experiments}, and the paper is concluded in \secref{sec:conclusions}.

\section{Preliminaries}\label{sec:models}
\subsection{Signals on the Sphere}
We consider square integrable complex functions $f(\theta,\phi)$ on the unit sphere $\untsph$, where $\theta\in[0,\pi]$ is the co-latitude and $\phi\in[0,2\pi)$ the longitude, with the inner product~\cite{Kennedy-book:2013}
\begin{align}\label{eqn:innprd}
\innerp{f}{h} \dfn \intsph f(\theta,\phi)\,\conj{h(\theta,\phi)}\,\sin\theta\,d\theta\,d\phi ,
\end{align}
where $\conj{(\cdot)}$ denotes complex conjugation and $\intsph=\inttheta\intphi$. The space of such functions with finite induced norm $\norm{f}\dfn\innerp{f}{f}^{1/2}$ is the Hilbert space $\lsph$.

\subsection{Harmonic Analysis on the Sphere}
The spherical harmonics $Y_\ell^m(\theta,\phi)$ of degree $\ell\ge0$ and order $|m|\le\ell$ form a complete orthonormal basis of $\lsph$ (see \appref{App:maths}), so that every $f\in\lsph$ has the expansion
\begin{equation}\label{Eq:f_expansion}
f(\theta,\phi)=\sum_{\ell=0}^{\infty}\sum_{m=-\ell}^{\ell}\shc{f}{\ell}{m}\,Y_\ell^m(\theta,\phi),\qquad
\shc{f}{\ell}{m}\dfn\innerp[\big]{f}{Y_\ell^m}.
\end{equation}
A signal is band-limited at $L$ if $\shc{f}{\ell}{m}=0$ for $\ell\ge L$; these signals form the $L^2$-dimensional subspace $\lsphL{L}$. We write $\plms_\ell^m(\theta)\dfn Y_\ell^m(\theta,0)$ for the scaled associated Legendre functions, so that $Y_\ell^m(\theta,\phi)=\plms_\ell^m(\theta)e^{im\phi}$, and $\plms_\ell^{-m}=(-1)^m\plms_\ell^m$ by \eqref{eq:yellm:conj}. Throughout, $x\dfn\cos\theta$ and $\mathcal{P}_n$ denotes the space of complex polynomials of degree at most $n$.

\subsection{Order-$m$ Profiles}
For $f\in\lsphL{L}$ and $|m|<L$ define the order-$m$ profile
\begin{align}\label{Eq:gm_integral}
G_m(\theta)\dfn\intphi f(\theta,\phi)\,e^{-im\phi}\,d\phi = 2\pi\sum_{\ell=|m|}^{L-1}\shc{f}{\ell}{m}\,\plms_\ell^m(\theta),
\end{align}
so that, by \eqref{Eq:f_expansion},
\begin{align}\label{Eq:f_from_G}
f(\theta,\phi)&=\frac{1}{2\pi}\sum_{m=-(L-1)}^{L-1}G_m(\theta)\,e^{im\phi},\\
\shc{f}{\ell}{m}&=\inttheta G_m(\theta)\,\plms_\ell^m(\theta)\,\sin\theta\,d\theta ,\label{Eq:coef_from_G}
\end{align}
which follows from $\shc{f}{\ell}{m}=\innerp{f}{Y_\ell^m}$ by integrating first over $\phi$. Inserting \eqref{Eq:gm_integral} and using $2\pi\inttheta(\plms_\ell^m)^2\sin\theta\,d\theta=\norm{Y_\ell^m}^2=1$ confirms the normalization.

\subsection{Optimal-Dimensionality Sampling Scheme}\label{sec:sampling_one}
\begin{figure}[t]
\centering
\IfFileExists{figures/fig02_sample_positions.pdf}{\includegraphics[width=\columnwidth]{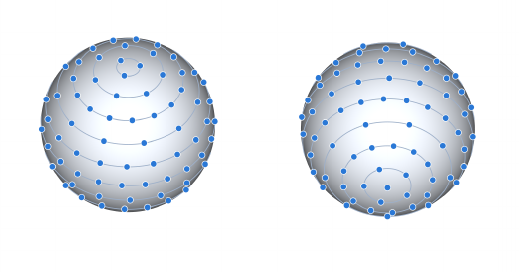}}{\placeholder{4.5cm}{Sample positions on the sphere for $L=12$ with the $\amp$-placement, viewed from the north pole and from the south pole (style of Fig.~2 of~\cite{Khalid:2014}).}}
\caption{Sample positions of the optimal-dimensionality sampling scheme for $L=12$ with the $\amp$-placement of \secref{sec:placement}, viewed from the north pole (left) and from the south pole (right), with $\nps=L^2=144$ samples; the single-point ring $k=0$ is the south pole $\theta_0=\pi$.}
\label{fig:sample_positions}
\end{figure}

Define the indexed vector $\bv\theta\dfn[\theta_0,\theta_1,\hdots,\theta_{L-1}]^T$ of $L$ co-latitudes and, for $|m|<L$, the sub-vector $\bv\theta^m\dfn[\theta_{|m|},\theta_{|m|+1},\hdots,\theta_{L-1}]^T$. The ring placed at $\theta_k$ carries $2k+1$ equally spaced points along $\phi$,
\begin{equation}
\bv\phi^k\dfn[\phi^k_0,\phi^k_1,\hdots,\phi^k_{2k}]^T,\qquad \phi^k_j\dfn\frac{2\pi j}{2k+1},
\end{equation}
so that the total number of samples is $\sum_{k=0}^{L-1}(2k+1)=L^2=\nps$. The co-latitudes are drawn from the equiangular set
\begin{equation}\label{Eq:equiangular_set}
\Theta=\{\vartheta_0,\vartheta_1,\hdots,\vartheta_{L-1}\},\qquad \vartheta_t\dfn\frac{\pi(2t+1)}{2L-1},
\end{equation}
and $\bv\theta$ is a particular ordering of $\Theta$ (each $\theta_k$ is one $\vartheta_t$), with the single-point ring at the pole $\theta_0=\pi$ so that $\sin\theta_k>0$ for all $k\ge1$. The \emph{alternating ordering} of~\cite{Khalid:2014} is the ordering in which $\theta_k$ is the $(k+1)$-th most distant element of $\Theta$ from the equator, so that the pole comes first and the two hemispheres alternate inward. Up to a reflection through the equator and a shift of the longitude origin by $\pi$ on alternate rings, this ordering is the ``1-3-5 distribution'' recommended in~\cite{Giacaglia:1972}. Both operations change only the signs of the rows and columns of the matrices $\Pm$ of \eqref{Eq:Pm}, so that the two schemes have the same condition numbers, which grow from about $10^4$ at $L=64$ to about $10^9$ at $L=128$~\cite{Khalid:2014}. We call the rings $k\ge|m|$ the \emph{resolving rings} of order $m$, since each carries at least $2|m|+1$ points, and the rings $k<|m|$ the \emph{small rings} of order $m$. The choice of the ordering is the subject of \secref{sec:placement}. Throughout, the \emph{$\kap$-placement} denotes the greedy placement of~\cite{Khalid:2014}, which minimizes the condition number $\kap_m$ of $\Pm$, and the \emph{$\amp$-placement} denotes the placement proposed in this paper, which minimizes the amplification constant $\amp_m$ of \secref{sec:error}; \figref{fig:sample_positions} shows the samples for $L=12$ with the $\amp$-placement.

For $|m|<L$ define the vectors and the matrix
\begin{align}\label{Eq:Pm}
\mathbf{g}_m&\dfn G_m(\bv\theta^m),\qquad
\mathbf{f}_m\dfn\big[\shc{f}{|m|}{m},\hdots,\shc{f}{L-1}{m}\big]^T,\nonumber\\
(\Pm)_{k,\ell}&\dfn 2\pi\,\plms_\ell^m(\theta_k),
\end{align}
with $k=|m|,\hdots,L-1$ and $\ell=|m|,\hdots,L-1$, so that $\mathbf{g}_m=\Pm\mathbf{f}_m$ by \eqref{Eq:gm_integral}. The forward SHT of~\cite{Khalid:2014}, which we refer to as the \emph{elimination transform}, computes $\mathbf{g}_m$ for $m=L-1,L-2,\hdots,0$ by FFTs along the resolving rings, solves $\mathbf{g}_m=\Pm\mathbf{f}_m$ (and the same for $-m$), evaluates
\begin{equation}\label{Eq:signal_subset}
\tilf_m(\theta,\phi)\dfn\frac{1}{2\pi}\Big(e^{im\phi}G_m(\theta)+e^{-im\phi}G_{-m}(\theta)\Big)
\end{equation}
on the small rings from the solved coefficients, and subtracts it there, so that the small rings become band-limited along $\phi$ at the next order. The transform developed below keeps this sweep and replaces the solve.

\section{Sampling Functions on Optimal-Dimensionality Samples}\label{sec:sampling_functions}
This section develops the Lagrange basis of $\lsphL{L}$ on the $L^2$ samples and shows that, on an iso-latitude scheme, it decomposes order by order into weighted one-dimensional Lagrange bases, which are the tool of the transform in \secref{sec:transform}.

\subsection{Sampling Functions and the Collocation Matrix}
Index the $L^2$ samples by $q=1,\hdots,\nps$ and write $\unit{x}_q$ for the $q$-th sample point. Define the collocation matrix $\mathbf{Y}\in\mathbb{C}^{\nps\times\nps}$ with entries $Y_\ell^m(\unit{x}_q)$, rows indexed by $q$ and columns by $(\ell,m)$. The point set is \emph{admissible} if $\mathbf{Y}$ is invertible. The \emph{sampling functions} of the point set are the $\nps$ signals
\begin{equation}\label{Eq:W_def}
W_q(\theta,\phi)\dfn\sum_{\ell=0}^{L-1}\sum_{m=-\ell}^{\ell}\big[\mathbf{Y}^{-1}\big]_{(\ell,m),q}\,Y_\ell^m(\theta,\phi),
\end{equation}
which satisfy $W_q(\unit{x}_s)=\delta_{qs}$.
Three properties follow directly from this definition. First, every $f\in\lsphL{L}$ can be written in terms of its samples as
\begin{equation}\label{Eq:interp_identity}
f(\theta,\phi)=\sum_{q=1}^{\nps}f(\unit{x}_q)\,W_q(\theta,\phi),
\end{equation}
since both sides lie in $\lsphL{L}$ and agree on $\nps$ points. Second, the sampling functions are orthonormal over the sample set. Third, $\sum_qW_q\equiv1$, since the constant function lies in $\lsphL{L}$. The sample vector and the coefficient vector are therefore two coordinate systems for the same signal, related by $\mathbf{Y}$, and the forward SHT is the application of $\mathbf{Y}^{-1}$ to the samples. The transforms of this paper and of~\cite{Khalid:2014} are two ways of applying $\mathbf{Y}^{-1}$ without forming it.

The sampling functions were introduced in this generality in~\cite{Giacaglia:1972}, by analogy with the cardinal functions of the Whittaker--Shannon sampling theorem. We note two of their properties that place the scheme among other point sets on the sphere. Let $K_L(\unit{x},\unit{y})=\sum_{\ell<L}\frac{2\ell+1}{4\pi}P_\ell(\unit{x}\cdot\unit{y})$ be the reproducing kernel of $\lsphL{L}$ and $\mathbf{K}$ its Gram matrix on the sample set, $\mathbf{K}_{qs}=K_L(\unit{x}_q,\unit{x}_s)=(\mathbf{Y}\mathbf{Y}^H)_{qs}$. Then $W_q(\unit{x})=\sum_s[\mathbf{K}^{-1}]_{qs}K_L(\unit{x},\unit{x}_s)$ and $\innerp{W_q}{W_s}=[\mathbf{K}^{-1}]_{qs}$. Hence the sampling functions are orthogonal in $\lsph$ exactly when $K_L$ vanishes between distinct samples, that is, when $P_{L-1}(\gamma)=P_L(\gamma)$ for every pairwise inner product $\gamma=\unit{x}_q\cdot\unit{x}_s$, $q\ne s$ (see the closed form of $K_L$ in \appref{App:maths}); in that case the sample set is an equal-weight quadrature of degree $2L-2$ with $L^2$ points, a tight spherical design~\cite{Delsarte:1977}, which on $\untsph$ exists, apart from the trivial case $L=1$, only for $L=2$, the regular tetrahedron~\cite{Bannai:2009}. The kernel form of the sampling functions and the condition $P_{L-1}(\gamma)=P_L(\gamma)$ are eqs.~(2.68) and~(2.71) of~\cite{Giacaglia:1972}, where the tetrahedron was found for $L=2$ and no configuration was found for $L=3$. Hence no placement of $L^2$ points, iso-latitude or not, can make the sampling functions orthogonal for $L\ge3$.

\subsection{Order-$m$ Structure on Iso-Latitude Rings}
Order-$m$ profiles are characterized by the following lemma.
\begin{lemma}[Order-$m$ profiles]\label{lem:space}
For $f\in\lsphL{L}$ and $|m|<L$, the profile $G_m$ of \eqref{Eq:gm_integral} belongs to
\begin{equation}
\Vm\dfn\big\{(\sin\theta)^{|m|}\,u(\cos\theta)\ :\ u\in\mathcal{P}_{L-1-|m|}\big\},
\end{equation}
and $\dim\Vm=L-|m|$.
\end{lemma}
\begin{IEEEproof}
By the Rodrigues formula, $\plms_\ell^m(\theta)=a_{\ell m}(1-x^2)^{|m|/2}\frac{d^{|m|}}{dx^{|m|}}P_\ell(x)=(\sin\theta)^{|m|}u_{\ell-|m|}(x)$ with $a_{\ell m}$ a constant and $u_{\ell-|m|}$ a polynomial of exact degree $\ell-|m|$. Hence $\{\plms_\ell^m\}_{\ell=|m|}^{L-1}$ is a basis of $\Vm$ and \eqref{Eq:gm_integral} places $G_m$ in it.
\end{IEEEproof}

A function in $\Vm$ is determined by its values at $L-|m|$ distinct co-latitudes away from the poles, and the resolving rings of order $m$ are exactly $L-|m|$ in number. This is the content of the next lemma, which defines the sampling functions of order $m$.
\begin{lemma}[Sampling functions of order $m$]\label{lem:lagrange}
Let $|m|<L$ and let $\theta_{|m|},\hdots,\theta_{L-1}$ be distinct, with $\sin\theta_k>0$ for $k\ge|m|$ when $m\ne0$. For $k\ge|m|$ define
\begin{equation}\label{Eq:Lambda_def}
\sfn{m}{k}(\theta)\dfn\Big(\frac{\sin\theta}{\sin\theta_k}\Big)^{|m|}\prod_{\substack{s=|m|\\ s\ne k}}^{L-1}\frac{\cos\theta-\cos\theta_s}{\cos\theta_k-\cos\theta_s}.
\end{equation}
Then $\sfn{m}{k}\in\Vm$, $\sfn{m}{k}(\theta_s)=\delta_{ks}$ for $k,s\ge|m|$, and for every $g\in\Vm$ and every $\theta\in[0,\pi]$
\begin{equation}\label{Eq:interp_m}
g(\theta)=\sum_{k=|m|}^{L-1}g(\theta_k)\,\sfn{m}{k}(\theta).
\end{equation}
\end{lemma}
\begin{IEEEproof}
The product in \eqref{Eq:Lambda_def} is the Lagrange polynomial $\lambda_k(x)\in\mathcal{P}_{L-1-|m|}$ on the $L-|m|$ distinct nodes $x_s\dfn\cos\theta_s$, so $\sfn{m}{k}=(\sin\theta)^{|m|}\lambda_k(x)/(\sin\theta_k)^{|m|}\in\Vm$ and $\sfn{m}{k}(\theta_s)=\lambda_k(x_s)=\delta_{ks}$. For $g=(\sin\theta)^{|m|}u$ with $u\in\mathcal{P}_{L-1-|m|}$, the right side of \eqref{Eq:interp_m} is $(\sin\theta)^{|m|}\sum_ku(x_k)\lambda_k(x)=(\sin\theta)^{|m|}u(x)$ by Lagrange interpolation. For $m=0$ no division by $\sin\theta_k$ occurs and the pole node is admissible.
\end{IEEEproof}

\begin{figure}[t]
\centering
\IfFileExists{figures/fig03_sampling_functions.pdf}{\includegraphics[width=\columnwidth]{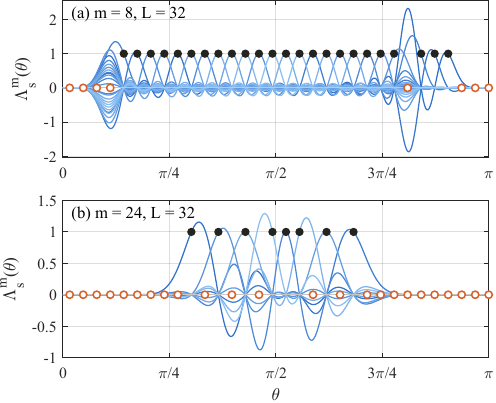}}{\placeholder{5.2cm}{Sampling functions $\sfn{m}{k}(\theta)$ of order $m=8$ for $L=32$ with the $\amp$-placement of \secref{sec:placement}: the $L-m=24$ functions on $[0,\pi]$, with the resolving co-latitudes $\theta_{8},\hdots,\theta_{31}$ marked by filled circles (each function equals one there and zero at the others) and the small-ring co-latitudes $\theta_0,\hdots,\theta_7$ marked by open circles, where the functions are evaluated by the transform. Second panel: the same for $m=24$, showing the $(\sin\theta)^{m}$ decay toward the poles.}}
\caption{Sampling functions $\sfn{m}{s}(\theta)$ of Lemma~\ref{lem:lagrange} for $L=32$ with the $\amp$-placement: (a) the $24$ functions of order $m=8$ and (b) the $8$ functions of order $m=24$, one curve for each resolving co-latitude $\theta_s$, $s\ge m$. Filled markers show the resolving co-latitudes, where each function equals one at its own node and zero at the others, and open markers show the co-latitudes of the small rings, where the transform evaluates the functions. Toward the poles the functions decay as $(\sin\theta)^{m}$, much faster for $m=24$ than for $m=8$. Between the nodes they can exceed one, by up to a factor of about $2.3$ in (a), next to the small ring that lies among the resolving rings near $\theta=0.8\pi$.}
\label{fig:sampling_functions}
\end{figure}

The sampling functions \eqref{Eq:Lambda_def} are eq.~(4.26) of~\cite{Giacaglia:1972}, where they were derived as the sampling functions of the space of associated Legendre functions of fixed order on the nodes $\theta_{|m|},\hdots,\theta_{L-1}$. The intuition is that of the cardinal series of Shannon restricted to a finite-dimensional space: the order-$m$ profile has $L-|m|$ degrees of freedom (the coefficients $\shc{f}{\ell}{m}$, $\ell=|m|,\hdots,L-1$), the resolving rings supply exactly $L-|m|$ values, and \eqref{Eq:interp_m} is the unique element of $\Vm$ through those values, with $\sfn{m}{k}$ playing the role of the sinc functions. In matrix terms, with $\Psm$ the $|m|\times(L-|m|)$ Legendre block on the small rings, $(\Psm)_{k,\ell}=2\pi\plms_\ell^m(\theta_k)$ for $k<|m|$,
\begin{equation}\label{Eq:Lambda_matrix}
\big[\Lm\big]_{k,s}\dfn\sfn{m}{s}(\theta_k)=\big[\Psm\,\Pm^{-1}\big]_{k,s},\quad k<|m|\le s ,
\end{equation}
so that the sampling functions of order $m$ are the closed form of the product of the two matrices that the elimination transform forms by solving. Throughout, $\sfn{m}{k}$ with a superscript denotes a sampling function, $\Lm$ in bold the matrix \eqref{Eq:Lambda_matrix} of its values, and $\amp_m$ with a subscript the amplification constant of \secref{sec:error}. The sampling functions are illustrated in \figref{fig:sampling_functions}.

\section{Proposed Optimal-Dimensionality Spherical Harmonic Transform}\label{sec:transform}
\subsection{Forward Transform}\label{sec:forward}
The proposed OSHT is summarized in \figref{fig:pipeline}. Let $r(\theta,\phi)$ denote a residual signal defined on the $\nps$ samples, initialized to $r=f$. The forward transform runs over the orders $m=L-1,L-2,\hdots,0$ and, at each order, treats $m$ and $-m$ together. Three steps are carried out.

\emph{(a) Resolve:} For every resolving ring $k\ge m$, compute
\begin{equation}\label{Eq:resolve}
G_{\pm m}(\theta_k)=\frac{2\pi}{2k+1}\sum_{j=0}^{2k}r(\theta_k,\phi_j^k)\,e^{\mp im\phi_j^k},
\end{equation}
which is the $\pm m$ bin of a $(2k+1)$-point FFT of the ring.

\emph{(b) Interpolate:} For every small ring $k<m$, compute
\begin{equation}\label{Eq:interpolate}
G_{\pm m}(\theta_k)=\sum_{s=m}^{L-1}G_{\pm m}(\theta_s)\,\sfn{m}{s}(\theta_k)
\end{equation}
from the values of step (a), using the sampling functions of order $m$ of Lemma~\ref{lem:lagrange}.

\emph{(c) Clean:} For every small ring $k<m$ and every point of that ring, update
\begin{equation}\label{Eq:clean}
r(\theta_k,\phi)\leftarrow r(\theta_k,\phi)-\frac{1}{2\pi}\Big(G_m(\theta_k)e^{im\phi}+G_{-m}(\theta_k)e^{-im\phi}\Big),
\end{equation}
which removes the orders $\pm m$ from the small rings. The sweep must proceed from the highest order to the lowest because of aliasing. A ring of $2k+1$ points folds every order $m'$ onto the bin $m'\bmod(2k+1)$, so that a high order contaminates the bins of lower orders, whereas a low order never contaminates the bin of a higher order. The only ring that can be read without any prior knowledge is therefore the largest ring, and only for the highest order. Before a smaller ring can be read at a lower order, the orders already resolved must be removed from it, and this is what step (c) does with the values provided by step (b). Steps (a) and (c) are those of~\cite{Khalid:2014}; step (b) replaces the solution of $\mathbf{g}_{\pm m}=\mathbf{P}_{\pm m}\mathbf{f}_{\pm m}$ and the evaluation of \eqref{Eq:signal_subset} from the solved coefficients. At the end of the sweep the profiles $G_m(\theta_k)$ are known for every $|m|<L$ and every ring $k$, those on the resolving rings from (a) and those on the small rings from (b), and no harmonic coefficient has yet been computed. The sweep is summarized in Procedure~\ref{Procedure:IE}.

\begin{algorithm}
\caption{Forward OSHT}\label{Procedure:IE}
\begin{algorithmic}[1]
\Require $f(\theta,\phi)$ on the $\nps$ samples; the ordering $\bv\theta$
\Procedure{Forward OSHT}{$f$}
\State $r\leftarrow f$
\For{$m=L-1,L-2,\hdots,0$}
\State compute $G_{\pm m}(\theta_k)$ for all $\theta_k\in\bv\theta^m$ by \eqref{Eq:resolve}
\NoNumber{(a $(2k+1)$-point FFT along each resolving ring)}
\State compute $G_{\pm m}(\theta_k)$ for all $\theta_k\in\bv\theta\setminus\bv\theta^m$
\NoNumber{by \eqref{Eq:interpolate}, with $\sfn{m}{s}(\theta_k)$ evaluated as in \secref{sec:numerics}}
\State update $r$ on all $\theta_k\in\bv\theta\setminus\bv\theta^m$ by \eqref{Eq:clean}
\EndFor
\State compute $\shc{f}{\ell}{m}$ from $\{G_m(\theta_k)\}$ by Procedure~\ref{Procedure:backend}
\State \textbf{return} $\shc{f}{\ell}{m}$
\EndProcedure
\end{algorithmic}
\end{algorithm}

\subsection{Back-Ends: From Profiles to Coefficients}\label{sec:backends}
Equation \eqref{Eq:coef_from_G} expresses each coefficient as an integral over $\theta$ of $G_m\plms_\ell^m$ against $\sin\theta$. By Lemma~\ref{lem:space} this product is $(1-x^2)^{|m|}$ times a polynomial of degree at most $(L-1-|m|)+(\ell-|m|)$ in $x$, hence a polynomial in $x$ of degree at most $2L-2$, and $\sin\theta\,d\theta=-dx$. Two exact evaluations follow.

\subsubsection{GL Back-End}\label{sec:gl_backend}
Let $x_n$, $n=1,\hdots,L$, be the $L$ zeros of the Legendre polynomial $P_L$ (\appref{App:maths}), which lie in $(-1,1)$, and let
\begin{equation}\label{Eq:GL_weights}
w_n\dfn\frac{2}{(1-x_n^2)\,\big[P_L'(x_n)\big]^2}>0
\end{equation}
be the corresponding weights, so that $\sum_{n=1}^{L}w_n\,u(x_n)=\int_{-1}^{1}u(x)\,dx$ for every polynomial $u$ of degree at most $2L-1$ (the $L$-point GL rule). Write $\theta_n^{\rm GL}\dfn\arccos x_n$. Evaluate the profile at the Gauss nodes by \eqref{Eq:interp_m}, $G_m(\theta_n^{\rm GL})=\sum_{s\ge|m|}G_m(\theta_s)\sfn{m}{s}(\theta_n^{\rm GL})$, and set
\begin{equation}\label{Eq:GL_backend}
\shc{f}{\ell}{m}=\sum_{n=1}^{L}w_n\,G_m(\theta_n^{\rm GL})\,\plms_\ell^m(\theta_n^{\rm GL}),\qquad |m|\le\ell<L .
\end{equation}
This back-end uses only the values on the resolving rings, and its weights are positive.

\subsubsection{Equiangular Back-End}\label{sec:eq_backend}
The ring co-latitudes $\Theta$ in \eqref{Eq:equiangular_set} are the grid of the equiangular sampling theorem~\cite{McEwen:2011}, whose forward algorithm takes as input exactly the profiles $G_m(\vartheta_t)$ on all $L$ rings. That algorithm can be understood as an exact quadrature in $\theta$ carried out in the Fourier domain of the torus $[0,2\pi)^2$. The optimal-dimensionality scheme supplies the profiles on only part of the torus grid, since the small rings carry fewer points, but the sweep reconstructs the missing values exactly, and the quadrature can then proceed as if $2L^2$ samples had been taken. After the sweep the profiles on all rings are available, and the algorithm applies without modification in four steps. First, each profile is extended periodically to $[0,2\pi)$ by $\tilde G_m(\vartheta_t)=G_m(\vartheta_t)$ for $t<L$ and $\tilde G_m(\vartheta_t)=(-1)^mG_m(\vartheta_{2L-2-t})$ for $t=L,\hdots,2L-2$, where $\vartheta_t=\pi(2t+1)/(2L-1)$ is now taken for $t=0,\hdots,2L-2$. Second, the extended profile is Fourier transformed in $\theta$, $F_{mm'}=\frac{1}{2\pi(2L-1)}\sum_{t=0}^{2L-2}\tilde G_m(\vartheta_t)e^{-im'\vartheta_t}$ for $|m'|<L$. Third, the quadrature is carried out implicitly by the convolution $H_{mm'}=2\pi\sum_{|m''|<L}F_{mm''}\,\varpi(m''-m')$, where $\varpi(m')\dfn\inttheta\sin\theta\,e^{im'\theta}\,d\theta$ equals $\pm i\pi/2$ for $m'=\pm1$, $2/(1-m'^2)$ for even $m'$ and zero otherwise, so that $H_{mm'}=\inttheta G_m(\theta)e^{-im'\theta}\sin\theta\,d\theta$. Finally, the coefficients are obtained as
\begin{equation}\label{Eq:MW_backend}
\shc{f}{\ell}{m}=i^{m}\sqrt{\frac{2\ell+1}{4\pi}}\sum_{m'=-(L-1)}^{L-1}\Delta^\ell_{m'm}\,\Delta^\ell_{m'0}\,H_{mm'},
\end{equation}
where $\Delta^\ell_{m'm}\dfn d^\ell_{m'm}(\pi/2)$ are the Wigner $d$-functions at $\pi/2$, computed by the recursion of Trapani and Navaza~\cite{Trapani:2006} (see \appref{App:maths}). Equation \eqref{Eq:MW_backend} follows from \eqref{Eq:coef_from_G} by inserting the Fourier expansion \eqref{Eq:dl_fourier} of $\plms_\ell^m(\theta)=\sqrt{(2\ell+1)/4\pi}\,d^\ell_{m0}(\theta)$. Both back-ends are collected in Procedure~\ref{Procedure:backend}. The GL back-end is simpler; the equiangular back-end uses all rings, needs no evaluation at Gauss nodes, and inherits the FFT structure and the existing implementations of the equiangular transform.

\begin{algorithm}
\caption{Back-Ends}\label{Procedure:backend}
\begin{algorithmic}[1]
\Require $G_m(\theta_k)$ for all $|m|<L$ and all $\theta_k\in\bv\theta$
\Procedure{GL back-end}{}
\For{$m=-(L-1),\hdots,L-1$}
\State evaluate $G_m(\theta_n^{\rm GL})$, $n=1,\hdots,L$, by \eqref{Eq:interp_m}
\State compute $\shc{f}{\ell}{m}$, $|m|\le\ell<L$, by \eqref{Eq:GL_backend}
\EndFor
\EndProcedure
\Procedure{Equiangular back-end}{}
\For{$m=-(L-1),\hdots,L-1$}
\State extend $G_m$ to $[0,2\pi)$; FFT in $\theta$ to $F_{mm'}$; convolve with $\varpi$ to $H_{mm'}$
\EndFor
\For{$\ell=0,\hdots,L-1$}
\State compute $\Delta^\ell$ by the Trapani--Navaza recursion; evaluate \eqref{Eq:MW_backend}
\EndFor
\EndProcedure
\end{algorithmic}
\end{algorithm}

\subsection{Exactness}
\begin{lemma}[Ring transform]\label{lem:dft}
If $h(\phi)=\sum_{|m'|\le M}b_{m'}e^{im'\phi}$ is sampled at $Q\ge2M+1$ equispaced points $\phi_j=2\pi j/Q$, then $\frac{2\pi}{Q}\sum_{j=0}^{Q-1}h(\phi_j)e^{-im\phi_j}=2\pi b_m$ for every $|m|\le M$.
\end{lemma}
\begin{IEEEproof}
$\frac1Q\sum_je^{i(m'-m)\phi_j}=1$ if $m'\equiv m\pmod Q$ and $0$ otherwise; $|m'-m|\le2M<Q$ forces $m'=m$.
\end{IEEEproof}

\begin{lemma}[Sweep invariant]\label{lem:inv}
In exact arithmetic, at the start of the step of order $m$ in Procedure~\ref{Procedure:IE}, every ring $k$ contains orders $|m'|\le\max(k,m)$ only, and the values produced by steps (a) and (b) at that step are exact.
\end{lemma}
\begin{IEEEproof}
By induction downward in $m$. At $m=L-1$ nothing has been subtracted and every ring contains orders $|m'|\le L-1$ only, so the claim holds by \eqref{Eq:f_from_G}. Assume it at step $m$. Every ring $k\ge m$ then contains orders $|m'|\le k$ only, and Lemma~\ref{lem:dft} with $Q=2k+1\ge2m+1$ and $M=k$ shows that step (a) returns the exact $G_{\pm m}(\theta_k)$ for $k\ge m$; Lemma~\ref{lem:lagrange} applied to $g=G_{\pm m}\in\Vm$ on the nodes $\theta_m,\hdots,\theta_{L-1}$ shows that step (b) returns the exact $G_{\pm m}(\theta_k)$ for $k<m$; step (c) then removes exactly the orders $\pm m$ from every ring $k<m$, which contains orders $|m'|\le m$ only, so that at the start of step $m-1$ every ring $k$ contains orders $|m'|\le\max(k,m-1)$ only.
\end{IEEEproof}

\begin{lemma}[Exact quadrature]\label{lem:gl}
For $f\in\lsphL{L}$ and $|m|\le\ell<L$, \eqref{Eq:GL_backend} holds exactly.
\end{lemma}
\begin{IEEEproof}
By \eqref{Eq:coef_from_G}, $\shc{f}{\ell}{m}=\inttheta G_m\plms_\ell^m\sin\theta\,d\theta=\int_{-1}^{1}u(x)\,dx$ with $u\in\mathcal{P}_{2L-2}$ by Lemma~\ref{lem:space}; the $L$-point GL rule \eqref{Eq:GL_weights} is exact to degree $2L-1$.
\end{IEEEproof}

\begin{theorem}[Exact recovery]\label{thm:exact}
Let $\bv\theta$ be an ordering of $\Theta$ with $\theta_0=\pi$. For every $f\in\lsphL{L}$, Procedure~\ref{Procedure:IE} followed by either back-end of Procedure~\ref{Procedure:backend} returns the exact spherical harmonic coefficients $\shc{f}{\ell}{m}$ in exact arithmetic. Consequently the restriction of $\lsphL{L}$ to the $\nps$ samples is injective, i.e., the optimal-dimensionality samples admit a sampling theorem for $\lsphL{L}$ in the collocation sense.
\end{theorem}
\begin{IEEEproof}
By Lemma~\ref{lem:inv}, after the sweep $G_m(\theta_k)$ is exact for all $|m|<L$ and all $k$. The GL back-end uses the resolving-ring values only, and Lemmas~\ref{lem:lagrange} and~\ref{lem:gl} give the exact coefficients. The equiangular back-end uses the values on all $L$ rings, which are the exact profiles $G_m(\vartheta_t)$ of a signal in $\lsphL{L}$ on its grid $\Theta$, and its forward algorithm is exact for such profiles~\cite{McEwen:2011}. Injectivity follows because the map from samples to coefficients just constructed is a left inverse of synthesis.
\end{IEEEproof}

\begin{remark}
The values on the small rings computed in step (b) are not used by the GL back-end, where their only role is to clean those rings so that step (a) is alias-free at the next order. The equiangular back-end uses them directly and needs no further interpolation. Theorem~\ref{thm:exact} also applies to the elimination transform of~\cite{Khalid:2014} with exact solves, since by \eqref{Eq:Lambda_matrix} the two transforms realize the same linear map at every order. The two transforms differ in how this map is computed: the proposed OSHT evaluates it in closed form through the sampling functions, and it obtains the coefficients by exact quadrature or by the equiangular transform rather than by matrix inversion.
\end{remark}

\subsection{Inverse Transform}\label{sec:inverse}
The inverse SHT is computed in two steps. First, for each $m$, the profile $G_m(\theta_k)=2\pi\sum_{\ell\ge|m|}\shc{f}{\ell}{m}\plms_\ell^m(\theta_k)$ is computed on all rings by the Legendre recurrence. Second, the signal $f(\theta_k,\phi)=\frac{1}{2\pi}\sum_mG_m(\theta_k)e^{im\phi}$ is evaluated on each ring by an inverse FFT of length $2k+1$, after the orders $|m|>k$ have been folded modulo $2k+1$. The folding is exact, since it reproduces the aliasing that the forward transform undoes. The cost of the inverse transform is $O(L^3)$.

\section{Numerical Realization and Computational Complexity}\label{sec:numerics}
\subsection{Stable Evaluation of the Sampling Functions}\label{sec:stable_eval}
The entries of $\Lm$ are moderate in size (\secref{sec:error}) but the two factors of \eqref{Eq:Lambda_def} are not: $(\sin\theta/\sin\theta_s)^{|m|}$ ranges over hundreds of orders of magnitude, and the Lagrange polynomial $\lambda_s(x)$ at a point outside the cluster of resolving nodes ranges over the inverse of that. Neither factor can therefore be formed on its own. We write, with $\omega_m(x)\dfn\prod_{s\ge|m|}(x-x_s)$ and the barycentric weights $\beta_s\dfn1/\prod_{s'\ge|m|,s'\ne s}(x_s-x_{s'})$, the first barycentric form $\lambda_s(x)=\omega_m(x)\beta_s/(x-x_s)$~\cite{Berrut:2004}, and evaluate
\begin{multline}\label{Eq:logLambda}
\log\abs{\sfn{m}{s}(\theta)}=|m|\big(\log\sin\theta-\log\sin\theta_s\big)\\
+\log\abs{\omega_m(x)}+\log\abs{\beta_s}-\log\abs{x-x_s},
\end{multline}
with the sign tracked separately, and exponentiate once. Every term is a sum of at most $L$ logarithms, so each entry is obtained with a relative error of order $L\epsilon_{\rm mach}$, up to the size of the logarithms involved, regardless of its magnitude, where $\epsilon_{\rm mach}$ is the unit roundoff, and no intermediate quantity can overflow or underflow. The first barycentric form is backward stable~\cite{Higham:2004}; the weighted ratio enters only inside \eqref{Eq:logLambda}. The small rings are never resolving nodes of the same order, and for every $L\le2048$ the Gauss nodes differ from all ring nodes by more than $0.18/L^2$ in $x$, so that $x-x_s\ne0$. A target exactly at a pole has $\sfn{m}{s}=0$ for $m\ne0$ and is treated as such. The unweighted alternative, which divides the profile by $(\sin\theta_s)^{|m|}$ to interpolate a plain polynomial, loses all accuracy, since the factors $(\sin\theta_s)^{-|m|}$ span many orders of magnitude.

The resolving node sets are nested: the set for order $|m|$ is the set for $|m|+1$ with $\theta_{|m|}$ added. The logarithms of the barycentric weights are therefore updated in $O(L-|m|)$ operations per order, $\log\abs{\beta_s}\leftarrow\log\abs{\beta_s}-\log\abs{x_s-x_{|m|}}$ for $s>|m|$ and $\log\abs{\beta_{|m|}}=-\sum_{s>|m|}\log\abs{x_{|m|}-x_s}$, so that the weights cost $O(L^2)$ for the whole sweep.

\subsection{Ring Transforms}\label{sec:ring}
The ring sums \eqref{Eq:resolve} are computed by FFTs of length $2k+1$, as in~\cite{Khalid:2014}. Since the samples on a ring have magnitude of order $L$ for coefficients of order one (the sum of $L^2$ terms of random sign), the root mean square rounding error of direct summation grows like $\sqrt{2k+1}$ times that magnitude, whereas that of the FFT grows like $\sqrt{\log(2k+1)}$ times it; in our experiments the FFT lowers the error of the whole transform by about a factor of two at $L\ge512$. One FFT per ring suffices for the entire sweep if the bins are kept and updated by the cleaning step, which changes only one bin of each ring at each order.

\subsection{Scaled Associated Legendre Functions}
The three-term recurrence in $\ell$ for fixed $m$, \eqref{Eq:recurrence} of \appref{App:maths}, is stable in the direction of increasing $\ell$ provided the seed $\plms_{m}^m(\theta)=(-1)^m\sqrt{\frac{(2m+1)!}{4\pi}}\frac{(\sin\theta)^{m}}{2^{m}m!}$, $m\ge0$, is computed accurately~\cite{Reinecke:2013}. For large $|m|$ and $\theta$ far from the equator the seed underflows double precision while $\plms_\ell^m(\theta)$ for larger $\ell$ is of order one, and a recurrence started from zero returns zero for the whole column. We therefore run the recurrence on a mantissa with a separate per-$\theta$ logarithmic scale, seeded by the logarithm of the seed computed through $\log\Gamma$, and rescale by $10^{-150}$ whenever the mantissa exceeds $10^{150}$; the stored value $\plms_\ell^m(\theta)$ is the product of mantissa and scale and underflows only when it is truly below $10^{-308}$. This is the adaptive rescaling used in the implementation of~\cite{Khalid:2014}, written out here in the form used by libsharp~\cite{Reinecke:2013}.

\subsection{Forbidden-Region Skipping}\label{sec:polar}
\begin{figure}[t]
\centering
\IfFileExists{figures/fig04_skipping.pdf}{\includegraphics[width=\columnwidth]{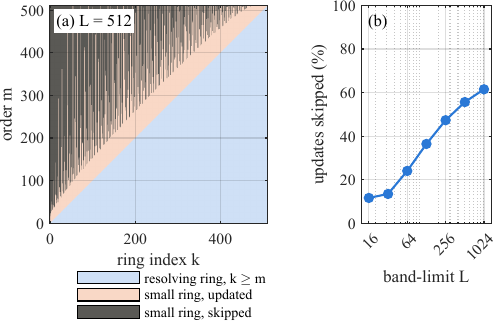}}{\placeholder{4.5cm}{Forbidden-region skipping: the $(m,k)$ plane for $L=512$ with the $\amp$-placement, cells shaded where the order-$m$ row of ring $k$ is below $\tau\varrho_L$ for $\tau=10^{-17}$ (skipped), the resolving region $k\ge|m|$ hatched; inset: fraction of small-ring updates skipped versus $L$ (56\% at $512$, 62\% at $1024$).}}
\caption{Forbidden-region skipping with $\tau=10^{-17}$ and the $\amp$-placement. (a) The pairs $(k,m)$ for $L=512$: resolving rings $k\ge m$ (light blue), small rings updated by steps (b) and (c) (light orange), and small rings whose order-$m$ row of Legendre functions lies below the threshold, for which these steps are skipped (dark). A small ring is skipped at the orders for which its co-latitude lies deep in the forbidden region, that is, at high orders and close to the poles. Since the placement does not order the rings by their distance from the equator, the boundary of the skipped region is irregular, and the light stripes inside it are small rings closer to the equator, which are updated up to higher orders. (b) Fraction of the small-ring updates that are skipped versus $L$, from about $12\%$ at $L=16$ to $62\%$ at $L=1024$.}
\label{fig:skipping}
\end{figure}

For $\theta$ beyond the turning point $\sin\theta\approx|m|/\ell$, the function $\plms_\ell^m(\theta)$ decays exponentially, and for a small ring near a pole the whole row $\{\plms_\ell^m(\theta_k)\}_{\ell<L}$ may lie below the double-precision resolution of the samples. Let $\tau$ be a relative threshold and $\varrho_L\dfn\sqrt{(2L-1)/4\pi}$ the scale of the largest $\plms$. We treat $\plms_\ell^m(\theta)$ as zero when $\abs{\plms_\ell^m(\theta)}<\tau\varrho_L$, and we skip steps (b) and (c) for a small ring, and the Gauss node in \eqref{Eq:GL_backend}, whose entire order-$m$ row is below the threshold. With $\tau=10^{-17}$ this removes 56\% of the small-ring updates at $L=512$ and 62\% at $L=1024$ with no change of the computed coefficients to three digits, since the skipped quantities are below what double precision registers in the sums they enter. The criterion must be the magnitude of the Legendre functions themselves; the factor $(\sin\theta_k/\sin\theta_s)^{|m|}$ alone is not a valid proxy, because it is paired in \eqref{Eq:Lambda_def} with a Lagrange polynomial that grows by the inverse amount, and a rule based on it destroys the transform. In a compiled implementation the decision can be taken from the turning point rather than by evaluating the row, and the saving then carries over to the Legendre recurrence, which can start at the turning point. The saving grows with $L$, since the decay beyond the turning point becomes steeper as the degree increases, so that a larger part of the forbidden region falls below the threshold. The skipped region and the saving as a function of $L$ are shown in \figref{fig:skipping}.

\subsection{Computational Complexity}\label{sec:complexity}
Step (a) costs one FFT of length $2k+1$ per ring, $O(L^2\log L)$ in total, when the bins are kept as in \secref{sec:ring}; recomputing the one bin needed on each resolving ring at each order costs $O(L^3)$ instead. Step (b) costs $\sum_m|m|(L-|m|)\approx L^3/3$ multiply-add operations over both signs of $m$, plus $O(L^2)$ for the nested barycentric weights and $\sum_m|m|(L-|m|)$ logarithms for \eqref{Eq:logLambda}. Step (c) costs $\sum_mm^2\approx2L^3/3$ multiply-add operations, also over both signs of $m$, when it is applied to the samples, and $O(L^2)$ when it is applied to the kept bins. The GL back-end costs $\sum_mL(L-|m|)\approx L^3$ evaluations of \eqref{Eq:logLambda} and the same number of Legendre values and multiply-add operations; the equiangular back-end costs $O(L^2\log L)$ in FFTs and $O(L^3)$ in the Wigner sums. The forward transform is therefore $O(L^3)$ with no linear solve and no precomputation beyond $\bv\theta$; the elimination transform of~\cite{Khalid:2014} has in addition $\sum_m(L-|m|)^3/3\approx L^4/6$ for the LU factorizations (both signs of $m$). Forbidden-region skipping halves the constants of steps (b) and (c). With $J$ passes of the multi-pass refinement of \secref{sec:multipass}, the cost is $J$ forward and $J-1$ inverse transforms. The comparison of the proposed scheme with the other schemes is summarized in \tabref{tab:complexity}.

The cost can be reduced further when the signal is real-valued. For real $f$ the profiles satisfy $G_{-m}=\conj{G_m}$ and the coefficients satisfy $\shc{f}{\ell}{-m}=(-1)^m\conj{\shc{f}{\ell}{m}}$, so that only the orders $m\ge0$ need to be processed and the cleaning step \eqref{Eq:clean} reduces to $r\leftarrow r-\frac1\pi\mathrm{Re}[G_m(\theta_k)e^{im\phi}]$.

\begin{table}[t]
\caption{Sampling schemes on the sphere for signals band-limited at $L$. Complexity is that of transforms with no precomputation other than the one listed; the fast algorithm of~\cite{Healy:2003} for the Driscoll--Healy scheme costs $O(L^2\log^2L)$ but requires $O(L^3)$ precomputation and storage.}
\label{tab:complexity}
\centering
\begin{tabular}{@{}lcccc@{}}
\toprule
Scheme & Samples & Complexity & Precomp. & Exact \\
\midrule
Driscoll--Healy~\cite{Driscoll:1994} & $\sim4L^2$ & $O(L^3)$ & none & yes \\
GL~\cite{Skukowsky:1986} & $\sim2L^2$ & $O(L^3)$ & nodes & yes \\
Equiangular~\cite{McEwen:2011} & $\sim2L^2$ & $O(L^3)$ & none & yes \\
Elimination~\cite{Khalid:2014} & $L^2$ & $O(L^4)$ & $\bv\theta$ & not claimed \\
This paper & $L^2$ & $O(L^3)$ & $\bv\theta$ & yes \\
\bottomrule
\end{tabular}
\end{table}

\section{Error Analysis, Placement of Rings and Multi-Pass Refinement}\label{sec:error}
\subsection{Amplification Constants}
In floating point, rounding errors enter Procedure~\ref{Procedure:IE} through the FFTs of step (a), the sums \eqref{Eq:interpolate} of step (b) and the subtractions of step (c), and the back-end adds the error of its own sums. Let $\xi_m$ bound the absolute error of the order-$m$ values on the resolving rings, which collects the errors of the FFTs and of the cleaning at earlier orders. The small-ring values then carry error at most $\amp_m\xi_m$ and the Gauss-node values at most $\amp_m^{\rm GL}\xi_m$, where
\begin{align}\label{Eq:amp}
\amp_m&\dfn\max_{k<|m|}\sum_{s=|m|}^{L-1}\abs{\sfn{m}{s}(\theta_k)}=\norm{\Psm\Pm^{-1}}_\infty,\\
\amp_m^{\rm GL}&\dfn\max_{1\le n\le L}\sum_{s=|m|}^{L-1}\abs{\sfn{m}{s}(\theta_n^{\rm GL})} .\label{Eq:ampGL}
\end{align}
The quadrature \eqref{Eq:GL_backend} amplifies by at most $\sum_nw_n\abs{\plms_\ell^m(\theta_n^{\rm GL})}\le\big(\sum_nw_n(\plms_\ell^m(\theta_n^{\rm GL}))^2\big)^{1/2}\big(\sum_nw_n\big)^{1/2}=(1/2\pi)^{1/2}\,2^{1/2}=\pi^{-1/2}$, using exactness of the rule for $(\plms_\ell^m)^2$, so the back-end contributes at most $\pi^{-1/2}\amp_m^{\rm GL}\xi_m$ to the coefficient error. No condition number is squared at any stage.

The constant $\amp_m$ has a simple interpretation. It is the Lebesgue constant of extrapolation from the resolving rings to the small rings in the weighted space $\Vm$, that is, the worst growth of a function of unit size on the resolving rings when it is continued to a small ring. For large $|m|$ the resolving rings cluster near the equator and the small rings lie near the poles, where the weight $(\sin\theta)^{|m|}$ suppresses the continuation. The dangerous case is a small ring at mid-latitude surrounded by resolving rings, where the continuation is an interpolation between widely spaced nodes, and the placement of \secref{sec:placement} is designed to avoid this case. In the elimination transform the corresponding constant is the condition number $\kap_m\dfn\kappa(\Pm)$, which bounds the relative error of the solve, and the two constants are related by $\amp_m\le\norm{\Psm}_\infty\norm{\Pm^{-1}}_\infty$. In our experiments the maximum of $\amp_m$ exceeds that of $\kap_m$ by a factor that grows with $L$, from about $1$ at $L=16$ to between $12$ and $18$ at $L=1024$, while the two transforms attain comparable accuracy. This indicates that $\amp_m$ is the relevant constant for both transforms: the error of the elimination transform is also governed by the map $\Psm\Pm^{-1}$, which that transform applies in two steps instead of one. In our experiments both constants are largest at small orders, near $m=50$, the constant $\amp_m^{\rm GL}$ stays close to $\amp_m$ at every order, and from $L=16$ to $L=1024$ the maximum of $\amp_m$ grows from about $7$ to $1940$ with the $\amp$-placement and to $3750$ with the $\kap$-placement.


\subsection{Error Propagation Through the Sweep}\label{sec:propagation}
That an error made at one order propagates to the orders processed after it was observed in~\cite{Khalid:2014}, where the error of the coefficients was shown to grow as the order decreases. We now describe the mechanism. The error that step (c) leaves on a small ring $k$ at order $m$ is a pure exponential $e^{\pm im\phi}$ of amplitude at most $\amp_m\xi_m/2\pi$. Since ring $k$ carries $2k+1$ points, this exponential aliases onto the single bin $\pm m\bmod(2k+1)$ of the ring, and it therefore corrupts exactly one lower order $m'$ on ring $k$. When the sweep reaches order $m'$, the corrupted value enters the interpolation of step (b), is amplified by about $\amp_{m'}$, and is passed on by step (c) to the rings $k'<|m'|$, where it corrupts a still lower order. The error therefore travels along a deterministic chain of orders, it is multiplied by an amplification constant at every link of the chain, and the number of links grows with $L$. This model explains two observations. First, the single-pass error grows much faster with $L$ than any single $\amp_m$ or $\kap_m$ once these constants are in the thousands; in our experiments it grows by a factor of about $300$ from $L=512$ to $L=1024$ with the $\kap$-placement and of about $80$ with the $\amp$-placement, and by a factor of about $9\times10^3$ from $L=1024$ to $L=2048$ with the $\amp$-placement, while $\max_m\amp_m$ grows only from $1940$ to about $2840$. Second, the error is of the same size for the elimination transform and the proposed OSHT, since the two share the chain.

\begin{figure}[t]
\centering
\IfFileExists{figures/fig05_error_vs_m.pdf}{\includegraphics[width=\columnwidth]{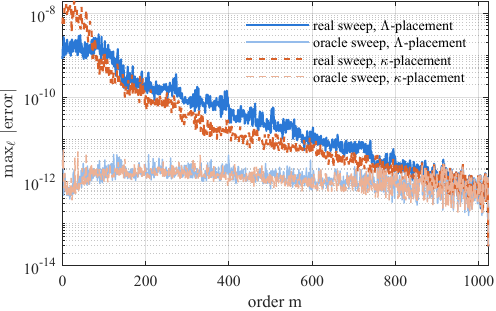}}{\placeholder{5cm}{Maximum coefficient error per order $m$ at $L=1024$, log scale: real sweep with the $\amp$-placement and with the $\kap$-placement, and oracle sweep with exact cleaning (flat at about $10^{-12}$).}}
\caption{Maximum coefficient error over $\ell$ and over $\pm m$ for each order $m$, for the real and the oracle sweeps of \secref{sec:propagation} with both placements, at $L=1024$ for one random test signal, one pass and the equiangular back-end. The real error reaches about $3\times10^{-9}$ with the $\amp$-placement and $2\times10^{-8}$ with the $\kap$-placement. The $\kap$-placement accumulates less error at intermediate orders but more at the smallest orders, which set the maximum error.}
\label{fig:error_vs_m}
\end{figure}

In order to verify the model, we carry out an experiment in which the propagation of error from one order to the next is switched off, which we refer to as the \emph{oracle} experiment. We generate a test signal band-limited at $L=1024$ with known coefficients, as described in \secref{sec:experiments}, and apply the forward transform to its samples twice. The first run is the transform of Procedure~\ref{Procedure:IE}, which we refer to as the real sweep. The second run, which we refer to as the oracle sweep, is identical to the first except in step (c): the profiles subtracted from the small rings are not the values interpolated in step (b) but the exact values $G_{\pm m}(\theta_k)$, computed from the known coefficients. In the oracle sweep the small rings are therefore cleaned without error, and the error in the coefficients of order $m$ is only the error made at order $m$ itself, by the FFTs of step (a), by the interpolation of step (b) and by the back-end. For each order $m$ we record the maximum error over the coefficients of orders $\pm m$ and plot it for both sweeps in \figref{fig:error_vs_m}. The error of the oracle sweep stays near $10^{-12}$ at every order. The error of the real sweep is about the same for the orders processed first, $m\ge900$, and then grows as the sweep proceeds toward $m=0$. With the $\amp$-placement its median value is about $6\times10^{-12}$ for $600\le m\le800$, $5\times10^{-11}$ for $300\le m\le500$ and $1.5\times10^{-9}$ for $m\le50$. Since the sweeps differ only in the cleaning step, the gap between the curves is the error carried down the chain, and the experiment confirms that the accuracy of the transform is set by the propagation of error through the sweep, not by the computation at any single order.

The model also indicates how the growth of the error can be limited. First, the placement of the rings sets the amplification constants $\amp_m$, and therefore the factor by which the error is multiplied at every link of the chain; this is the subject of \secref{sec:placement}. Second, every error made by a local step enters the chain, so each local step must be computed as accurately as possible; this is the purpose of the FFTs, the log-domain evaluation of the sampling functions and the scaled Legendre recurrence of \secref{sec:numerics}. Third, the error that remains after one pass is removed by the multi-pass refinement of \secref{sec:multipass}.

\subsection{Placement of Rings by Amplification Minimization}\label{sec:placement}
Since $\amp_m$ is the constant that governs the elimination, we place the rings to minimize it, which gives the $\amp$-placement. The construction follows the greedy procedure of~\cite{Khalid:2014}. The ring of $2L-1$ points is placed at the element of $\Theta$ nearest the equator, and for $m=L-2,L-3,\hdots,0$ the ring of $2m+1$ points is placed at the unused element $c\in\Theta$ that minimizes
\begin{equation}\label{Eq:greedy_objective}
\max_{\theta'\in U\setminus\{c\}}\ \sum_{s\in\{c\}\cup\{\theta_{m+1},\hdots,\theta_{L-1}\}}\abs{\sfn{m}{s}(\theta')},
\end{equation}
where $U$ is the set of unused elements of $\Theta$; these are exactly the small rings of order $m$ in the final ordering, so \eqref{Eq:greedy_objective} is $\amp_m$ for the candidate. Each candidate costs $O(m(L-m))$ through \eqref{Eq:logLambda}, so the search costs $O(L^4)$ in all and is vectorized; it is a precomputation performed once per $L$, and the result is stored as the vector $\bv\theta$ (Procedure~\ref{Procedure:placement}). In our experiments the $\amp$-placement lowers $\max_m\amp_m$ by a factor $7$ at $L=512$ and also lowers $\max_m\kap_m$ (from $209$ to $61$), so that the single-pass error of the transform improves by a factor of $4$ to $13$ for $256\le L\le1024$, and that of the elimination transform by a similar factor; at $L=2048$ the improvement is about $5\times10^3$ (\secref{sec:experiments}). Both greedy placements keep the largest rings at the co-latitudes nearest the equator, alternating between the two hemispheres as the alternating ordering does, and depart from that ordering mainly for the rings of intermediate size, which they move to mid-latitudes. For $L\le64$ the two placements give worst cases within a factor of two of each other, and the same worst case at $L=64$; the gain of the $\amp$-placement grows with $L$. The worst case is attained at a small order ($m\le51$ for every $L\le1024$), where the greedy procedure has the fewest co-latitudes left to choose from, and a look-ahead search or an elimination-style search~\cite{Nafees:2017} over \eqref{Eq:greedy_objective} is the natural next step.

\begin{algorithm}
\caption{Placement of Rings by Amplification Minimization}\label{Procedure:placement}
\begin{algorithmic}[1]
\Require $L$
\Procedure{Placement}{$L$}
\State $U\leftarrow\Theta$ of \eqref{Eq:equiangular_set}; $\theta_{L-1}\leftarrow$ element of $U$ nearest $\pi/2$; $U\leftarrow U\setminus\{\theta_{L-1}\}$
\For{$m=L-2,L-3,\hdots,0$}
\State for each $c\in U$ evaluate \eqref{Eq:greedy_objective} by \eqref{Eq:logLambda}
\State $\theta_m\leftarrow\arg\min_c$; $U\leftarrow U\setminus\{\theta_m\}$
\EndFor
\State \textbf{return} $\bv\theta=[\theta_0,\hdots,\theta_{L-1}]^T$
\EndProcedure
\end{algorithmic}
\end{algorithm}

For completeness we note that the condition-number placement of~\cite{Khalid:2014}, the $\kap$-placement, can be computed in $O(L^4)$ rather than $O(L^5)$: at step $m$ every candidate matrix is a fixed block $\mathbf{B}$ (the rows of the rings already placed) plus one candidate row $\bv{v}^T$, so $\Pm^H\Pm=\mathbf{B}^H\mathbf{B}+\conj{\bv v}\bv{v}^T$; diagonalizing $\mathbf{B}^H\mathbf{B}$ once per order, at a cost of $O((L-m)^3)$, and finding the extreme eigenvalues of the rank-one update for every candidate by bisection on the secular equation, at a cost of $O((L-m)^2)$ per candidate, gives all condition numbers of the order in $O((L-m)^3+m(L-m)^2)$ operations, and the whole placement in $O(L^4)$. This implementation reproduces, ring for ring, the stored orderings of the implementation of~\cite{Khalid:2014} at $L=256$ and $L=721$, and we used it to compute the $\kap$-placement up to $L=2048$.

\subsection{Multi-Pass Refinement}\label{sec:multipass}
We adopt the multi-pass SHT of~\cite{Nafees:2017} without change. Let $\hat f_1$ denote the coefficients obtained from one pass of the forward transform. We synthesize $\hat f_1$ on the samples by the inverse transform, form the residual $r_1=f-\hat f_1$ on the samples, apply the forward transform to the residual, and set $\hat f_2=\hat f_1+\widehat{r_1}$. The procedure is repeated until the maximum of the residual on the samples stops decreasing. This is iterative refinement in the sense of linear algebra. The transform computes $\mathbf{Y}^{-1}$ with a relative error $\eta$ that depends on the sweep but not on the signal, and the residual is itself a signal in $\lsphL{L}$ whose magnitude is the error of the previous pass. Transforming the residual reproduces this error with the same relative accuracy, so that each pass multiplies the error by $\eta$, until the error reaches a floor. The floor is set by the rounding of the samples, about $\epsilon_{\rm mach}\max|f|$, amplified by $\mathbf{Y}^{-1}$. Since the synthesis is exact up to rounding, the refinement converges to the coefficients of the band-limited signal that interpolates the samples, which is $f$ itself. The amplification by $\mathbf{Y}^{-1}$ grows with $L$ on $L^2$ samples. The condition number of $\mathbf{Y}$, scaled by $(4\pi/\nps)^{1/2}$ so that an ideal equal-weight quadrature would give $1$, is $82$ at $L=64$ for both placements, whereas the $\ngl\times L^2$ and $\nmw\times L^2$ collocation matrices of the GL and equiangular schemes have condition numbers $5.2$ and $10.7$. This is the price of having no redundancy in the samples; a fuller analysis of the stability and geometry of the sampling scheme is beyond the scope of this paper.

\section{Numerical Experiments}\label{sec:experiments}

\begin{figure*}[t]
\centering
\subfloat[]{\IfFileExists{figures/fig06a_emean.pdf}{\includegraphics[width=0.32\textwidth]{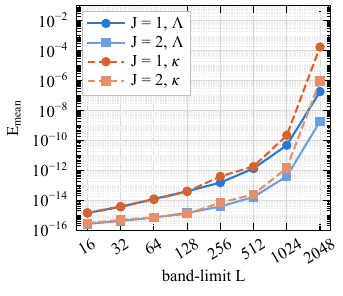}}{\parbox{0.32\textwidth}{\placeholder{4.3cm}{$E_{\rm mean}$ versus $L$ after $J=1,2$ passes}}}\label{fig:multipass_emean}}
\hfill
\subfloat[]{\IfFileExists{figures/fig06b_emax.pdf}{\includegraphics[width=0.32\textwidth]{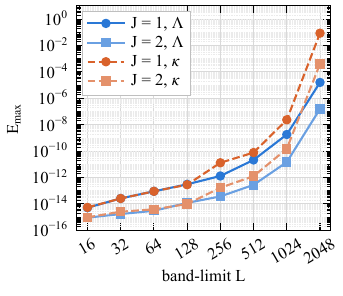}}{\parbox{0.32\textwidth}{\placeholder{4.3cm}{$E_{\rm max}$ versus $L$ after $J=1,2$ passes}}}\label{fig:multipass_emax}}
\hfill
\subfloat[]{\IfFileExists{figures/fig06c_snr.pdf}{\includegraphics[width=0.32\textwidth]{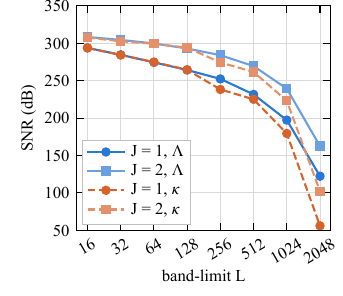}}{\parbox{0.32\textwidth}{\placeholder{4.3cm}{$\mathrm{SNR}$ versus $L$ after $J=1,2$ passes}}}\label{fig:multipass_snr}}
\caption{Accuracy of the proposed OSHT with the equiangular back-end in a single pass ($J=1$) and after a second pass of the multi-pass refinement ($J=2$): (a) mean error $E_{\rm mean}$, (b) maximum error $E_{\rm max}$ and (c) SNR versus band-limit, averaged over 5 test signals, for the $\amp$-placement and the $\kap$-placement, denoted $\amp$ (solid lines) and $\kap$ (dashed lines) in the legends.}
\label{fig:multipass}
\end{figure*}

\begin{figure}[t]
\centering
\IfFileExists{figures/fig07_timing.pdf}{\includegraphics{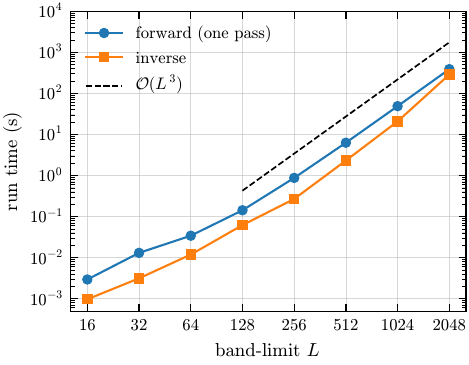}}{\placeholder{4.5cm}{Computation time versus $L$ (log-log) of one forward pass and of the inverse transform, $16\le L\le2048$, with an $O(L^3)$ reference line.}}
\caption{Computation time versus band-limit of one pass of the proposed forward OSHT, from the samples to the coefficients, with the equiangular back-end and the $\amp$-placement, and of the inverse transform, with an $O(L^3)$ reference line.}
\label{fig:timing}
\end{figure}

The accuracy of the forward transform is measured on band-limited test signals whose coefficients are known. For each test signal we draw the coefficients $\shc{f_{\rm t}}{\ell}{m}$, $0\le\ell<L$, $|m|\le\ell$, with real and imaginary parts independent and uniformly distributed in $[-1,1]$, and synthesize $f_{\rm t}$ on the $\nps$ samples by the inverse transform of \secref{sec:inverse}, which is exact up to rounding. The forward transform is then applied to these samples and returns the coefficients $\shc{f_{\rm r}}{\ell}{m}$. Since $f_{\rm t}\in\lsphL{L}$, an exact transform would return $\shc{f_{\rm r}}{\ell}{m}=\shc{f_{\rm t}}{\ell}{m}$ (Theorem~\ref{thm:exact}), so every difference between the two sets of coefficients is numerical error of the transform. We record
\begin{align}\label{Eq:exp1:errors}
E_{\rm max}&\dfn\max\abs{\shc{f_{\rm t}}{\ell}{m}-\shc{f_{\rm r}}{\ell}{m}},\\
E_{\rm mean}&\dfn\frac{1}{L^2}\sum_{\ell,m}\abs{\shc{f_{\rm t}}{\ell}{m}-\shc{f_{\rm r}}{\ell}{m}},
\end{align}
together with the signal-to-noise ratio (SNR),
\begin{equation}\label{Eq:snr}
\mathrm{SNR}\dfn-10\log_{10}\frac{\sum_{\ell,m}\abs{\shc{f_{\rm t}}{\ell}{m}-\shc{f_{\rm r}}{\ell}{m}}^2}{\sum_{\ell,m}\abs{\shc{f_{\rm t}}{\ell}{m}}^2}\ \text{dB},
\end{equation}
averaged over 5 signals, for $L=16,32,\hdots,2048$. The maximum error is the quantity used in~\cite{Khalid:2014,Nafees:2017} and is the most severe, since it is set by the worst coefficient; the SNR puts the energy of the error in relation to the energy of the signal and is the figure of merit commonly used in applications. The accuracy is reported for the proposed OSHT with the equiangular back-end, with the $\kap$-placement and the $\amp$-placement and with one and two passes, and with three passes up to $L=1024$. Computation times are measured for one pass of the forward transform, from the samples to the coefficients, with the equiangular back-end and the $\amp$-placement, and for the inverse transform, as the smallest of up to five repeated runs. The results of Figs.~\ref{fig:multipass} and~\ref{fig:timing} are computed in double precision with the Python package \texttt{osht} (Python 3.13.2, NumPy 2.5.4, SciPy 1.18.1) on a laptop computer with an Intel Core i7-1365U processor (10 cores, 12 threads, up to 5.2 GHz) and 32 GB of memory, running Windows 11. The accuracy is reported in \figref{fig:multipass} and the computation time in \figref{fig:timing}; Figs.~\ref{fig:sample_positions} to~\ref{fig:error_vs_m} in the earlier sections are also produced by these experiments.

The accuracy results of \figref{fig:multipass} show the following. (i) For $L\le128$ the two placements give the same error, and for $256\le L\le1024$ the $\amp$-placement lowers the single-pass maximum error by a factor of $4$ to $13$, from $1.3\times10^{-11}$ to $1.3\times10^{-12}$ at $L=256$ and from $2.3\times10^{-8}$ to $1.8\times10^{-9}$ at $L=1024$. At $L=2048$ it lowers the error from $8.6\times10^{-2}$ to $1.6\times10^{-5}$, and the SNR of a single pass from $56$~dB to $122$~dB. (ii) A second pass reduces the maximum error at $L=1024$ to $1.5\times10^{-11}$ with the $\amp$-placement and $1.5\times10^{-10}$ with the $\kap$-placement, and a third pass brings no significant further reduction, since the floor of \secref{sec:multipass} has been reached. At $L=2048$ the second pass reduces the maximum error to $1.6\times10^{-7}$ with the $\amp$-placement and $4.1\times10^{-4}$ with the $\kap$-placement. (iii) The time of one forward pass grows by factors of $7.2$, $7.8$ and $8.0$ over the last three doublings of $L$, in agreement with the complexity $O(L^3)$, and is $49$~s at $L=1024$ and $390$~s at $L=2048$; the inverse transform takes $21$~s and $290$~s (\figref{fig:timing}).

\section{Conclusions}\label{sec:conclusions}
We have developed a spherical harmonic transform for the optimal-dimensionality sampling scheme that is theoretically exact, solves no linear system, requires no precomputation beyond the placement of the rings, and has complexity $O(L^3)$. The transform rests on the sampling functions of order $m$, a weighted Lagrange basis on the rings that resolve order $m$, through which the profile of the signal on every other ring follows in closed form. The harmonic coefficients are then obtained by GL quadrature or by the equiangular transform, both of which are exact. We have identified the amplification constant that governs the numerical error of the elimination, used it to place the rings, explained the growth of the single-pass error with band-limit by a propagation model, and shown that the multi-pass refinement removes what remains in two passes. We have also placed the scheme in the context of the sampling functions of~\cite{Giacaglia:1972}, which contain the point set and the closed-form inverse but not the transform, and of spherical designs, which explain why orthogonal sampling functions on $L^2$ points cannot exist. The open problems are a bound on the amplification constants and on the length of the propagation chain, a placement that goes beyond the greedy at small orders, and the behavior of the transform at band-limits beyond those considered here, where the single-pass error is expected to keep growing and the accuracy will depend increasingly on the multi-pass refinement.

\appendices
\section{Mathematical Background}\label{App:maths}
The spherical harmonic of degree $\ell\ge0$ and order $|m|\le\ell$ is~\cite{Kennedy-book:2013,Sakurai:1994}
\begin{equation}\label{Eq:Sph_harmonics}
Y_\ell^m(\theta,\phi)=N_\ell^m\,P_\ell^m(\cos\theta)\,e^{im\phi},\ \ N_\ell^m\dfn\sqrt{\tfrac{2\ell+1}{4\pi}\tfrac{(\ell-m)!}{(\ell+m)!}},
\end{equation}
with $P_\ell^m$ the associated Legendre function, $P_\ell^m(x)=\frac{(-1)^m}{2^\ell\ell!}(1-x^2)^{m/2}\frac{d^{\ell+m}}{dx^{\ell+m}}(x^2-1)^\ell$ for $0\le m\le\ell$ and $P_\ell^{-m}=(-1)^m\frac{(\ell-m)!}{(\ell+m)!}P_\ell^m$, so that $\innerp{Y_\ell^m}{Y_{\ell'}^{m'}}=\delta_{\ell\ell'}\delta_{mm'}$ and
\begin{equation}\label{eq:yellm:conj}
\conj{Y_\ell^m(\theta,\phi)}=(-1)^mY_\ell^{-m}(\theta,\phi).
\end{equation}
Here $P_\ell\dfn P_\ell^0$ is the Legendre polynomial of degree $\ell$. The scaled functions $\plms_\ell^m(\theta)=N_\ell^mP_\ell^m(\cos\theta)$ of \secref{sec:models} satisfy, for fixed $m\ge0$ and $\ell\ge m$, the three-term recurrence
\begin{equation}\label{Eq:recurrence}
\plms_{\ell+1}^m=\mu_\ell^m\cos\theta\,\plms_\ell^m-\nu_\ell^m\,\plms_{\ell-1}^m,
\end{equation}
with
\begin{align}
\mu_\ell^m&=\sqrt{\frac{(2\ell+1)(2\ell+3)}{(\ell+1-m)(\ell+1+m)}},\nonumber\\
\nu_\ell^m&=\sqrt{\frac{(2\ell+3)(\ell-m)(\ell+m)}{(2\ell-1)(\ell+1-m)(\ell+1+m)}}
\end{align}
(so that $\nu_m^m=0$), started from $\plms_{m-1}^m=0$ and the seed $\plms_m^m$ given in \secref{sec:numerics}. The addition theorem gives the reproducing kernel of $\lsphL{L}$: $\sum_{\ell<L}\sum_{|m|\le\ell}Y_\ell^m(\unit{x})\conj{Y_\ell^m(\unit{y})}=\sum_{\ell<L}\frac{2\ell+1}{4\pi}P_\ell(\unit{x}\cdot\unit{y})=\frac{L}{4\pi}\frac{P_{L-1}(\gamma)-P_L(\gamma)}{1-\gamma}$ with $\gamma=\unit{x}\cdot\unit{y}$, by the Christoffel--Darboux formula.

The Wigner $d$-functions $d^\ell_{m'm}(\beta)$, $|m'|,|m|\le\ell$, are the real matrix elements of the rotation by the angle $\beta$ about the $y$-axis in the degree-$\ell$ representation, in the convention of~\cite{Sakurai:1994,McEwen:2011}, and satisfy $Y_\ell^m(\theta,\phi)=\sqrt{(2\ell+1)/4\pi}\,d^\ell_{m0}(\theta)e^{im\phi}$, that is, $\plms_\ell^m(\theta)=\sqrt{(2\ell+1)/4\pi}\,d^\ell_{m0}(\theta)$. They admit the Fourier expansion~\cite{Trapani:2006,McEwen:2011}
\begin{equation}\label{Eq:dl_fourier}
d^\ell_{m'm}(\beta)=i^{m-m'}\sum_{m''=-\ell}^{\ell}\Delta^\ell_{m''m'}\,\Delta^\ell_{m''m}\,e^{im''\beta},
\end{equation}
with $\Delta^\ell_{m'm}\dfn d^\ell_{m'm}(\pi/2)$, which is what turns the integral \eqref{Eq:coef_from_G} into the sums of \eqref{Eq:MW_backend}. The planes $\Delta^\ell$, $\ell=0,1,\hdots,L-1$, are computed by the recursion of Trapani and Navaza~\cite{Trapani:2006} at a cost of $O(\ell^2)$ per plane, and the recursion is stable up to at least $L=2048$~\cite{McEwen:2011}.


\end{document}